\documentstyle[12pt]{article}

\input psfig.tex 

\newcommand{\beq}{\begin{equation}}
\newcommand{\eeq}{\end{equation}}

\newcommand{\into}{\rightarrow}

\newtheorem{conjecture}{Conjecture}

\begin{document}

\pagestyle{plain}
\pagenumbering{arabic}
\setcounter{page}{1}

\title{On the Asymptotics of Takeuchi Numbers}

\author{Thomas Prellberg\thanks{email: {\tt tprell@phy.syr.edu}}
\thanks{Address after 1 June 2000: Institut f\"ur Theoretische Physik, Abteilung C, 
Technische Universit\"at Clausthal, D-38678 Clausthal-Zellerfeld, Germany}\\
Physics Department, Syracuse University\\
Syracuse, NY 13244, USA}

\date{May 1, 2000} 
 
\maketitle 
 
\begin{abstract}  
I present an asymptotic formula for the Takeuchi numbers $T_n$.
In particular, I give compelling numerical evidence and present a
heuristic argument showing that 
$$
T_n\sim C_T\,B_n\exp{1\over2}{W(n)}^2
$$
as $n$ tends to infinity, where $B_n$ are the Bell numbers, 
$W(n)$ is Lambert's $W$ function, and $C_T=2.239\ldots$ is
a constant. Moreover, I show that the method presented here can be generalized 
to derive conjectures for related problems. \end{abstract}

\newpage

\section{Introduction}

In a paper entitled ``Textbook Examples of Recursion,''
Donald E.\ Knuth discusses recurrence equations related to the
properties of recursive programs \cite{knuth91}, among them 
Takeuchi's function \cite{takeuchi78,takeuchi79}
\beq
t(x,y,z)=\mbox{\bf if $x\leq y$ then $y$ else 
$t(t(x-1,y,z),t(y-1,z,x),t(z-1,x,y))$}\;.
\eeq
Let $T(x,y,z)$ denote the number of times the {\bf else} clause is
invoked when $t(x,y,z)$ is evaluated recursively.\footnote{Note that
it is the recursive evaluation of $t(x,y,z)$ rather than the actual value of $t(x,y,z)$ 
that is of interest. See Knuth's paper for an explicit expression of $t(x,y,z)$.}
For non-negative integers $n$, the Takeuchi numbers $T_n$ are 
defined as $T_n=T(n,0,n+1)$.
The first few values of $T_n$ for $n=0,1,2,\ldots$ are
\beq
0,\;1,\;4,\;14,\;53,\;223,\;1034,\;5221,\;28437,\;165859,\;\ldots\;. 
\eeq
Knuth gives the recurrence
\beq
\label{takrec}
T_{n+1}=\sum_{k=0}^n\left\{{n+k\choose n}-{n+k\choose n+1}\right\}T_{n-k}
+\sum_{k=1}^{n+1}{2k\choose k}{1\over k+1}\quad\mbox{, $n\geq 0$,}
\eeq
and deduces a functional equation for the generating function 
$T(z)=\sum_{n=0}^\infty T_nz^n$:
\beq
\label{takfunc}
T(z)={C(z)-1\over1-z}+{z(2-C(z))\over\sqrt{1-4z}}T(zC(z))\;,
\eeq
where 
\beq
C(z)={1\over2z}(1-\sqrt{1-4z})=\sum_{n=0}^\infty {2n\choose n}{1\over n+1}
\eeq
is the generating function for the Catalan numbers
$C_n={2n\choose n}{1\over n+1}$.
Lastly, he gives asymptotically valid bounds for $T_n$,
\beq
\label{bounds}
e^{n\log n-n\log\log n-n}<T_n<e^{n\log n-n+\log n}\qquad\mbox
{for all sufficiently large $n$,}
\eeq
and poses obtaining further information about the
asymptotic properties of $T_n$ as an open problem. 

In this paper I give arguments leading to two conjectures about the
asymptotic behaviour of the Takeuchi numbers $T_n$. In Section 2, I 
present an explicit asymptotic formula for $T_n$ which improves upon the bounds 
(\ref{bounds}) based on numerical evidence and a heuristic argument.
For this, I briefly discuss the related asymptotic behavior of the Bell numbers 
and give an argument based on a numerical observation which leads directly to an 
explicit asymptotic formula for $T_n$ as $n$ tends to infinity. The
formula, as described in Conjecture 1, is exact up to $O((\log n/n)^2)$ and
contains a constant $C_T$ which is numerically determined to 25 significant digits. 
Section 3 presents a heuristic analytic argument which gives 
the asymptotic behavior up to $o(1)$ and enables one to identify the 
constant $C_T$ in terms of an explicit expression, as stated in Conjecture 2.  
In the final section I show that the method developed in Section 3 can give insight into 
the asymptotic behavior of a larger class of problems.

I conclude this introduction by briefly discussing the structure of the 
recurrence (\ref{takrec}) and the related functional equation 
(\ref{takfunc}). It is clear from the asymptotic bounds (\ref{bounds})
for $T_n$, that the generating function $T(z)$, defined as a formal power 
series, does not converge, and is therefore at best only an asymptotic 
expansion to an actual solution of the functional equation. This is 
also evident from the structure of the functional equation. This 
structure becomes clearer upon a change of variables,
\beq
T(z)={1\over z}T(z-z^2)-{1\over(1-z)(1-z+z^2)}\;,
\eeq
where one sees directly that the transformation involved is $g(z)=z-z^2$,
which is only marginally contracting at its fixed point $z=0$. While functional
equations with a transformation $g(z)$ that has an expanding or contracting fixed point 
($|g'(0)|\neq1$) are very well understood \cite{kuczma90}, it is precisely the fact that
$|g'(0)|=1$ which is at the root of the underlying difficulty of the problem discussed in this paper.

\section{Numerical Observations}

Our starting point is Knuth's observation \cite{knuth91} that for 
$n>0$ the Bell numbers $B_n$ are a lower bound to $T_n$. Here, $B_n$ is 
defined as
\beq
\label{bellrec}
B_{n+1}=\sum_{k=0}^n{n\choose k}B_{n-k}\;,\qquad B_0=1\;.
\eeq
The asymptotics of $B_n$ is discussed in great detail by de Bruijn 
\cite{debruijn61}. In \cite{moser55} one finds a systematic way of 
generating higher order terms in the asymptotic expansion by means 
of a contour integral representation using the well-known fact that
\beq
\label{bellfunc}
\sum_{n=0}^\infty B_n{z^n\over n!}=\exp(e^z-1)\;.
\eeq
Alternatively, one can also expand the right-hand side of (\ref{bellfunc})
in $z$ to get
\beq\label{bellsum}
B_n={1\over e}\sum_{m=0}^\infty{m^n\over m!}\;,
\eeq
which can then be evaluated asymptotically by use of the Euler-MacLaurin
formula. Either way, in the course of the computation of the asymptotics
of $B_n$ it turns out that a convenient asymptotic scale is given in
terms of $W(n)$ rather than in terms of $n$, 
where  $W(x)$ is Lambert's $W$ function, which is defined as the real
solution of
\beq
W(x)\exp W(x)=x\;.
\eeq 
The sum (\ref{bellsum}) is dominated by terms around $m=e^{W(n)}$, and one can easily calculate
that, written in terms of $w=W(n)$, the Bell numbers $B_n$ behave asymptotically 
as
\begin{eqnarray}
\label{bellasy}
\log B_n&=&
e^w(w^2-w+1)-{1\over2}\log(1+w)-1-{w(2w^2+7w+10)\over24(1+w)^3}e^{-w}\\
&&-{w(2w^4+12w^3+29w^2+40w+36)\over48(1+w)^6}e^{-2w}+O(e^{-3w})\nonumber\;,
\end{eqnarray}
and it is straightforward to calculate additional terms. It is more conventional
to state this formula with the exponentials $e^{kw}$ replaced by $(n/w)^k$, but this
obscures the fact that the asymptotic expansion is obtained in terms of $w$ rather than $n$.

Having rather explicit control over this lower bound, it is natural to
now try to compare $B_n$ and $T_n$ more closely. There 
is a principal difficulty coming from the fact that the asymptotic 
scale presumably also involves $w=W(n)$, which grows more slowly than $\log n$.
(As an example, $W(1000)\approx5.2496$ and 
$W(10000)\approx7.2318$.) Thus, one would expect that a direct numerical
investigation of $T_n/B_n$ is not very insightful, due to the presence 
of slowly varying correction terms of unknown form. 

\begin{figure}[th]
\label{figure1}
\centerline{\psfig{file=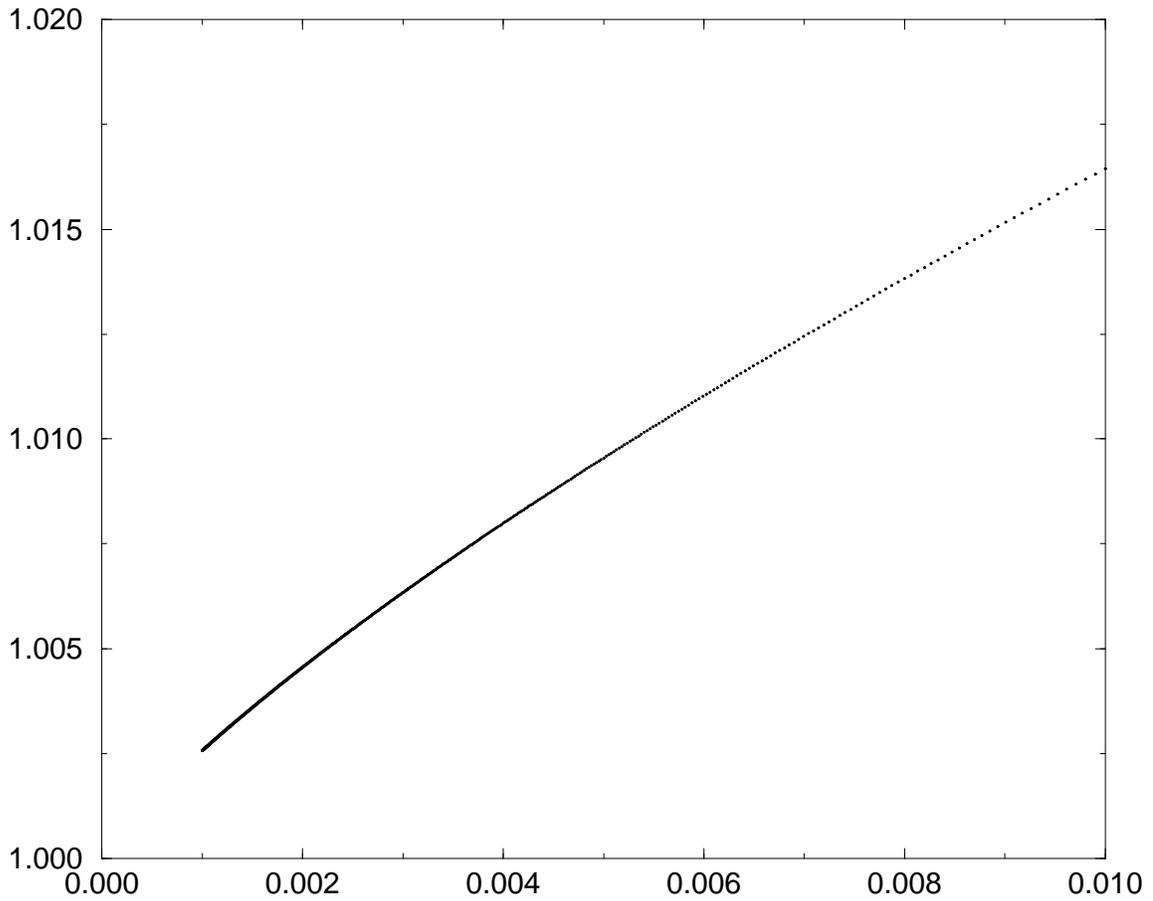,height=5in,angle=-90}}
\caption{${T_{n+1}/T_n}-{B_n/B_{n-1}}$ plotted versus $1/n$ for $n\leq1000$}
\end{figure}

\begin{figure}[th]
\label{figure2}
\centerline{\psfig{file=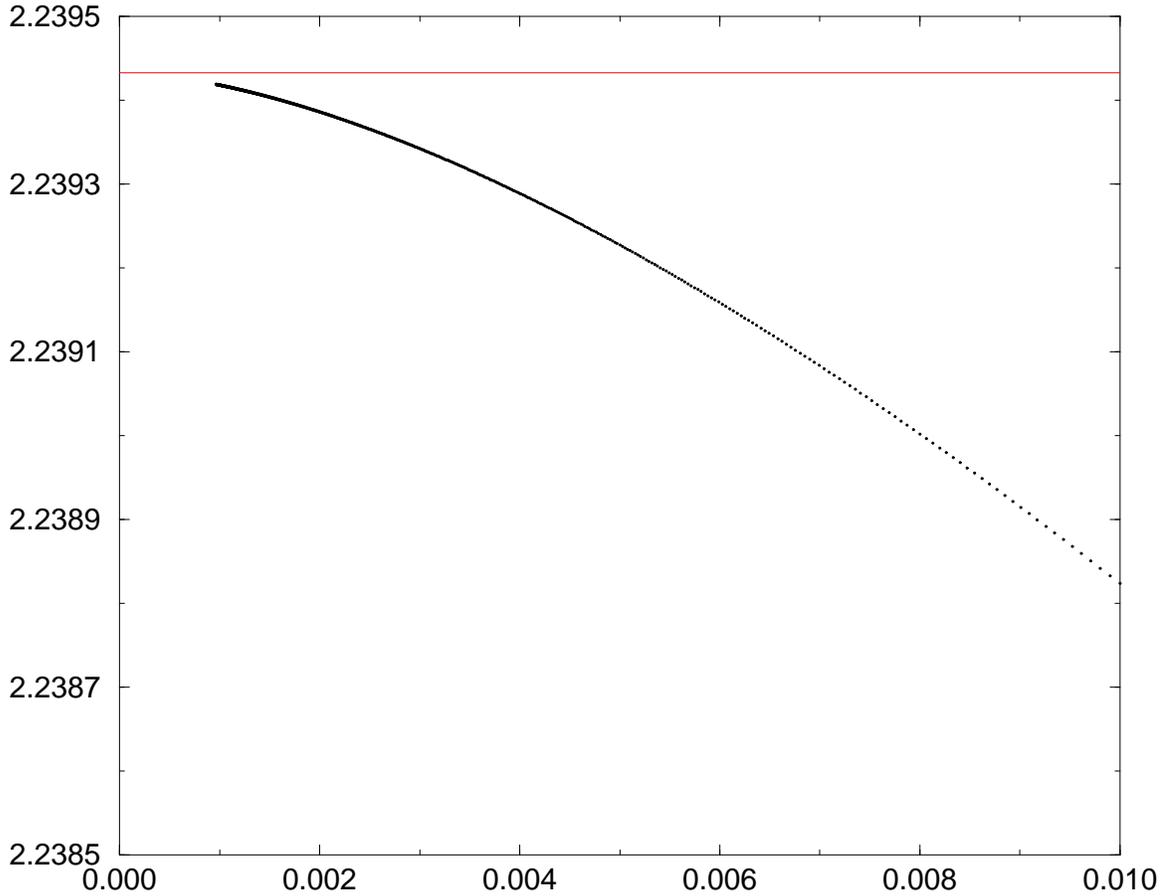,height=5in,angle=-90}}
\caption{${T_{n+1}/[B_n\exp({1\over2}w^2+w)]}$ plotted versus $1/n$ for $n\leq1000$. The horizontal line is
at $C_T=2.2394331040\ldots$.}
\end{figure}

However, as Figure 1 shows, if one compares the growth rates $T_n/T_{n-1}$ and $B_n/B_{n-1}$ 
instead, one is led to observe the surprisingly simple relationship
\beq
\label{numobs}
\lim_{n\into\infty}\left({T_{n+1}\over T_n}-{B_n\over B_{n-1}}\right)=1\;.
\eeq
In fact, the left hand side approaches $1$ rather quickly,
\beq
{B_n\over B_{n-1}}+1\leq{T_{n+1}\over T_n}
\leq{B_n\over B_{n-1}}+1+O(e^{-w})\;.
\eeq
This (unproven) numerical observation leads to a straightforward
derivation of an asymptotic formula. From (\ref{bellasy}) it follows 
easily that $B_{n-1}/B_n=e^{-w}+O(e^{-2w})$. Now one takes logarithms 
and sums up successively, from whence it follows that
\beq
\log T_{n+1}=\log B_n+{1\over2}w^2+w+O(1)\;.
\eeq

Now that I have guessed the leading asymptotic form, I can again resort to
numerical work to try to improve upon it. In fact, numerically it appears 
that the convergence is even better than expected due to a chance 
cancellation of higher order correction terms. Figure 2 shows that
\beq
T_{n+1}=C_T\,B_n\exp({1\over2}w^2+w+O(e^{-2w}))\;,
\eeq
and from the first $1000$ series terms I am able to deduce 
by iterative application of standard series extrapolation methods 
that 
\beq
C_T=2.23943\;31040\;05260\;73175\;4785\;(1)\;.
\eeq

Using the known asymptotic form of the Bell numbers, I can now give 
an explicit asymptotic expression for $T_n$ in terms of $w=W(n)$ alone,
as stated in the following conjecture.

\begin{conjecture}
As $n$ tends to infinity, one has
\begin{eqnarray}
\log T_n&=&e^w(w^2-w+1)+{1\over2}w^2-{1\over2}\log(1+w)+\\
&&+\log C_T-1-{w(26w^2+67w+46)\over24(1+w)^3}e^{-w}+O(e^{-2w})\nonumber\;.
\end{eqnarray}
Here $w=W(n)$ and the constant $C_T$ is some positive real number.
\end{conjecture}

Numerically, $\log C_T-1=-0.19377\;72447\;31916\;75890\;1157\;(1)$. Of course it would 
be desirable to find an analytic expression for this number. In the next 
section I shall present a heuristic argument giving such an expression. 

Dropping one correction term and comparing with the asymptotic 
expression (\ref{bellasy}) for the Bell numbers, Conjecture 1 implies the nice formula given 
in the abstract,
\beq
\label{takbell}
T_n\sim C_T\,B_n\exp{1\over2}{W(n)}^2\;.
\eeq

\section{Analytic Results}

In view of the previous section it seems promising to exploit the
apparent affinity between Takeuchi numbers $T_n$ and Bell numbers $B_n$. 
Given a recurrence of the general form
\beq
\label{rec}
a_n=\sum_{k=1}^nc_{n,k}a_{n-k}+b_n\;,
\eeq
I choose to write
\beq
\label{ansatz}
a_n={1\over e}\sum_{m=0}^\infty{1\over m!}f_{m,n}\quad\mbox{and}\quad b_n={1\over e}\sum_{m=0}^\infty{1\over m!}b_n\;.
\eeq
Inserting (\ref{ansatz}) into the recurrence (\ref{rec}) and shifting the 
summation index by one, I next equate terms to get
\beq
f_{m,n}=m\sum_{k=1}^nc_{n,k}f_{m-1,n-k}+b_n\;,
\eeq
This Ansatz might seem less arbitrary when considering that 
in the case of Bell numbers it reduces to $f_{m,n}=m^n$. In general, 
one observes that $f_{m,n}$ must be a polynomial in $m$ of at most degree $n$,
which I write as
\beq
\label{eq1}
f_{m,n}=m^n\sum_{k=0}^nd_{n,k}m^{-k}\;.
\eeq
If one further requires the coefficients $c_{n,k}$ in the recurrence
to be polynomials of degree $k$ in $n$, it follows that $d_{n,k}$ 
are polynomials of degree $k$ in $n$, which I write as
\beq
\label{eq2}
d_{n,k}=n^k\sum_{l=0}^kr_{k,l}n^{-l}\;.
\eeq
Combining equations (\ref{eq1}) and (\ref{eq2}) gives
\beq
m^{-n}f_{m,n}=\sum_{k=0}^n\sum_{l=0}^kr_{k,l}n^{k-l}m^{-k}
=\sum_{l=0}^nm^{-l}\sum_{k=0}^{n-l}r_{l+k,l}(n/m)^k\;.
\eeq
In order to get an idea about the asymptotic behavior of this double sum, 
I now replace the quotient $n/m$ by a new variable $v$ and consider the 
formal limit of taking the summation bounds to infinity, leading to
\beq
s_m(v)=\sum_{l=0}^\infty m^{-l}r_l(v)
\qquad\mbox{with}\qquad 
r_l(v)=\sum_{k=0}^\infty r_{l+k,l}v^k\;.
\eeq
Applying this method to Takeuchi numbers, one now inserts
$c_{n,k}=\left\{{n+k-2\choose n-1}-{n+k-2\choose n}\right\}$
and $b_n=\sum_{k=1}^n{2k\choose k}{1\over k+1}$.
With this choice, $r_0(v)$ is trivially zero as $c_{n,k}$ are polynomials 
in $n$ of degree $k-1$, and one gets a rather interesting result for $l\geq1$.
In fact,
\begin{eqnarray}
r_1(v)&=&e^{{1\over2}v^2+v}\\
r_2(v)&=&e^{{1\over2}v^2}\left(2e^{2v}-{1\over2}(v^3+v^2+4v+2)e^v\right)\\
r_3(v)&=&e^{{1\over2}v^2}\left(-{1\over8}e^{3v}-(v^3+3v^2+7v+6)e^{2v}+\right.\\
&&\left.+{1\over24}(3v^6+6v^5+47v^4+52v^3+144v^2+74v+51)e^v\right)\nonumber\\
r_4(v)&=&e^{{1\over2}v^2}\left(-{347\over108}e^{4v}
+{1\over16}(v^3+5v^2+12v+12)e^{3v}+\right.\\
&&+{1\over12}(3v^6+18v^5+89v^4+226v^3+411v^2+406v+195)e^{2v}-\nonumber\\
&&-{1\over432}(9v^9+27v^8+315v^7+603v^6+3024v^5+\nonumber\\
&&\left.+3384v^4+8757v^3+4707v^2+5484v+772)e^{v}\right)\nonumber\;.
\end{eqnarray}
From this, I conjecture
\beq
r_l(v)=e^{{1\over2}v^2}\left(p_{l,0}(v)e^{lv}+p_{l,1}(v)e^{(l-1)v}+\ldots
+p_{l,l-1}(v)e^{v}\right)\;,
\eeq
where $p_{l,k}(v)$ are polynomials in $v$ of degree $3k$. (This pattern
has been verified for $l\leq 8$.) For $v$ large, this conjecture implies
\beq
r_l(v)\sim\lambda_le^{{1\over2}v^2+lv}\;,
\eeq
and with a little effort one can compute the next values of 
$\lambda_l=p_{l,0}$. One gets
\beq
\lambda_0=0\;,\quad
\lambda_1=1\;,\quad
\lambda_2=2\;,\quad
\lambda_3=-{1\over8}\;,\quad
\lambda_4=-{347\over108}\;,\quad
\lambda_5={28201\over3456}\;,
\nonumber
\eeq
\beq
\lambda_6=-{3172987\over216000}\;,\quad
\lambda_7={822813607\over93312000}\;,\quad
\lambda_8={2183235065857\over16003008000}\;,\quad\ldots\nonumber\;,
\eeq
from whence I am led to conjecture that one can write 
$\lambda_n=\mu_n/[(n-1)!]^3$ where $\mu_n$ is integer.
(Unfortunately, I did not find a way to compute the coefficients 
$\lambda_l$ in a closed form!) One would expect from this behavior that
\beq
s_m(v)\sim e^{{1\over2}v^2}h(e^v/m)\;,
\eeq
where $h(x)$ is given in some sense by
\beq
h(x)=\sum_{l=0}^\infty\lambda_lx^l\;.
\eeq
I caution here that the series may be divergent and just valid as an 
asymptotic expansion. 

Keeping in mind that the evidence for the existence of $h(x)$ is rather 
sketchy, I nevertheless proceed under the assumption that for $n>>m>>1$ 
one can write 
\beq
f_{m,n}\sim m^ne^{{1\over2}(n/m)^2}h(e^{n/m}/m)\;.
\eeq
This now enables a heuristic computation of the Takeuchi numbers $T_n$.
I approximate
\begin{eqnarray}
T_n&=&{1\over e}\sum_{m=0}^\infty{1\over m!}f_{m,n}
\sim{1\over e}\sum_{m\approx m_{max}(n)}{1\over m!}f_{m,n}\\
&\sim&{1\over e}\sum_{m\approx m_{max}(n)}
{m^n\over m!}e^{{1\over2}(n/m)^2}h(e^{n/m}/m)\;,
\end{eqnarray}
where in the last step it is assumed that $n>>m_{max}(n)>>1$. This sum
is indeed dominated around $m_{max}(n)\sim e^w$, where the argument of $h$ 
simplifies to $1$. A careful asymptotic analysis of
\beq
\hat T_n={1\over e}\sum_{m\approx m_{max}(n)}
{m^n\over m!}e^{{1\over2}(n/m)^2}h(e^{n/m}/m)
\eeq
gives
\begin{eqnarray}
\label{that}
\log\hat T_n&=&e^w(w^2-w+1)+{1\over2}w^2-{1\over2}\log(1+w)+h_0-1\\
&&+{w(12w^5+24w^4+36w^3+58w^2+29w-10)\over24(w+1)^3}e^{-w}
\nonumber\\
&&+{(w+1)(h_1^2+h_2)+(2w^2+w+2)h_1\over2}e^{-w}+O(e^{-2w})\;,
\nonumber
\end{eqnarray}
where one has expanded $h(x)$ around $x=1$ as 
$\log h(x)=h_0+h_1(x-1)+h_2(x-1)^2/2+O((x-1)^3)$.

As long as the corrections made on passing from $T_n$ to $\hat T_n$
are small enough, it follows easily from this that asymptotically
\beq
T_n\sim B_ne^{{1\over2}w^2}h(1)\;,
\eeq
and one can identify the constant $C_T$ from equation (\ref{takbell}) 
with $h(1)$. Provided
the series expansion of $h(x)=\sum_{k=0}^\infty\lambda_kx^k$ converges
at $x=1$, I can thus conjecture an explicit expression for the constant 
$C_T$, which in principle is computable. 

\begin{conjecture}
The constant $C_T$ in Conjecture $1$ is given by
\beq
C_T=h(1)=\sum_{k=0}^\infty\lambda_k\;.
\eeq
\end{conjecture}

While the approximation of $T_n$ by $\hat T_n$ may be correct 
up to $O(e^{-w})$, no choice of $h(x)$ can match the next term in 
(\ref{that}) with the expansion of $T_n$. Thus, one also gets an 
indication of the size of the error made.

It seems that a careful asymptotic evaluation of the $f_{m,n}$ 
promises to be a suitable way of providing rigorous proof for the 
asymptotics of the Takeuchi numbers. Of course one could also try to
find a direct proof of our numerically observed equation (\ref{numobs}).

\section{A Generalization}

In the derivation of the functional equation (\ref{takfunc}) for the 
Takeuchi numbers $T_n$, it is crucial that
\beq
\sum_{k=0}^\infty{n+2k\choose k}z^k=C(z)^k/\sqrt{1-4z}\;,
\eeq
as this identity allows the explicit summation of the terms in the recurrence (\ref{takrec}).
The identity used is a special case of the following nice identity
\beq
\label{ident}
\sum_{k=0}^\infty{n+(\lambda+1)k\choose k}z^k=
\left\{\sum_{k=0}^\infty{(\lambda+1)k\choose k}z^k\right\}
\left\{\sum_{k=0}^\infty{(\lambda+1)k\choose k}
{z^k\over1+\lambda k}\right\}^n\;.
\eeq
This identity can be proved by inserting $z=y/(1+y)^{\lambda+1}$, which
after expanding leads to
\beq
\sum_{k=0}^\infty{(\lambda+1)k\choose k}{z^k\over1+\lambda k}=1+y
\eeq
and
\beq
\sum_{k=0}^\infty{n+(\lambda+1)k\choose k}z^k
={(1+y)^{n+1}\over1-\lambda y}\;.
\eeq
I use this now as a motivation for the study of the family of
recursions (with parameter $\lambda$)
\beq
\label{family}
A_{n+1}=\sum_{k=0}^n{n+\lambda k\choose k}A_{n-k}\;,\qquad A_0=1\;.
\eeq
Due to equation (\ref{ident}) one is able to derive a functional equation 
for the corresponding generating function $A(z)=\sum_{n=0}^\infty A_nz^n$:
\beq
A(z)=1+z{1+y\over1-\lambda y}A(z(1+y))\;,\quad z=y/(1+y)^{\lambda+1}\;.
\eeq
For $\lambda=0$ one recovers the recursion for the Bell numbers, and for
$\lambda=1$ one has something which is at least ``morally'' related to
the Takeuchi numbers.

Inserting the Ansatz (\ref{ansatz}) into (\ref{family}), one can easily 
repeat the analysis of the previous section. The result is now
\beq
A_n\sim B_n\exp\lambda\left\{{1\over2}{W(n)}^2+W(n)+d(\lambda)\right\}
\eeq
for any fixed value of $\lambda$. Again, one have an identification of the 
kind $d(\lambda)=h_\lambda(1)$, where the first terms
in the series expansion of $h_\lambda(x)$ are
\begin{eqnarray}
h_\lambda(x)&=&{1\over2}(\lambda-1)x-{1\over24}(2\lambda^2+18\lambda-5)x^2
-{1\over216}(33\lambda^3+90\lambda^2-329\lambda+54)x^3\nonumber\\
&-&{1\over960}(52\lambda^4-520\lambda^3+4240\lambda-502)x^4
+O(x^5)\;,
\end{eqnarray}
and one sees that the $k$th coefficient is a polynomial in $\lambda$ of 
degree $k$ (this has been verified up to $k=7$). I caution again that
convergence of this series expansion is an open question.

Finally, one can establish numerically the next term in the asymptotic
expansion of $A_n$. For any fixed value of $\lambda$, one finds
\beq
\log A_n=\log B_n+
\lambda\left({w^2\over2}+w+d(\lambda)-{\lambda+1\over2}e^{-w}\right)
+O(e^{-2w})\;.
\eeq
Indeed, this result even seems to hold for {\it complex} values 
of $\lambda$.

I conclude with remarking that even though Takeuchi's function has been
labelled a ``Textbook Example,'' it provides an exciting open question
for asympotic analysis.

\section*{Acknowledgements}

I thank Philippe Flajolet for bringing this problem to my attention.


\begin{thebibliography}{Takeuchi 1979}

\bibitem[DeBruijn 1961]{debruijn61} 
N. G. de Bruijn. 
{\it Asymptotic Methods in Analysis},
North Holland, Amsterdam, 1961.

\bibitem[Knuth 1991]{knuth91} 
Donald E. Knuth. 
Textbook Examples of Recursion. 
In {\it Artificial Intelligence and Theory of Computation}, 
Academic Press, London, pages 207-229, 1991. 

\bibitem[Kuczma 1990]{kuczma90}
M. Kuczma, B. Choczewski, and R. Ger,
{\it Iterative Functional Equations},
Encyclopedia of Mathematics and its Applications,
Cambridge University Press, 1990.

\bibitem[Moser 1955]{moser55}
Leo Moser and Max Wyman.
An Asymptotic Formula for the Bell Numbers.
In {\it Transactions of the Royal Society of Canada}
{\bf 49}, 49-53, 1955.

\bibitem[Takeuchi 1978]{takeuchi78} 
Ikuo Takeuchi. 
On a Recursive Function That Does Almost Recursion Only. 
Memorandum, Musahino Electrical Communication Laboratory, 
Nippon Telephone and Telegraph Co., Tokyo, 1978.

\bibitem[Takeuchi 1979]{takeuchi79} 
Ikuo Takeuchi. 
Dai-Ni-Kai LISP Kontesuto [On the Second LISP Contest]. 
In {\it J\=oh\=o Shori} 
{\bf 20}, pages 192-199, 1979.

\end{thebibliography}
\end{document}